\documentclass[smallextended]{svjour3}
\smartqed
\usepackage{theorem,amsmath,latexsym,amssymb}
\newcommand{\pf}{{\em Proof\/: }}

\newcommand{\Z}{\mathbb{Z}}
\newcommand{\F}{\mathbb{F}}
\newcommand{\supp}{{\rm supp}}
\newcommand{\soc}{{\rm soc}}
\newcommand{\rad}{{\rm rad}}

\newcommand{\res}{{\rm Res}}

\newcommand{\sh}{{\rm Sho}}

\newcommand{\xa}{\alpha}

\newcounter{abccnt}
\newenvironment{abc}{\begin{list}{\bf(\alph{abccnt})}{\usecounter{abccnt}
\labelwidth4ex \labelsep1ex \leftmargin6ex
\parsep3pt \itemsep1pt \topsep3pt}}{\end{list}}

\begin{document}

\title{New Bounds for Codes over Finite Frobenius Rings\thanks{Research supported Science Foundation Ireland Grant 08/RFP/MTH1181,Research supported by the Claude
Shannon Institute, Science Foundation Ireland Grant 06/MI/006} }

\author{Eimear Byrne \and Marcus Greferath \and Axel Kohnert \and  Vitaly Skachek}

\institute{E. Byrne, \at School of Mathematical Sciences
University College Dublin, Belfield, Dublin 4, Ireland\\
\email{ebyrne@ucd.ie}
\and
M. Greferath  \at School of Mathematical Sciences
University College Dublin, Belfield, Dublin 4, Ireland\\
\email{marcus.greferath@ucd.ie}
\and
A. Kohnert \at Mathematics Department, University of Bayreuth,
D-95440 Bayreuth, Germany\\
\email{axel.kohnert@uni-bayreuth.de}
\and
Vitaly Skachek
\at School of Mathematical Sciences
University College Dublin, Belfield, Dublin 4, Ireland\\
\email{vitaly.skackek@ucd.ie}
}

%\author{Eimear Byrne \\ \inst{School of Mathematical Sciences\\ 
%University College Dublin \\ Belfield, Dublin 4 \\ Ireland}
%\footnote{Research supported Science Foundation Ireland Grant 08/RFP/MTH1181}\\
%\and Marcus Greferath  \\ \inst{School of Mathematical Sciences\\ 
%University College Dublin \\ Belfield, Dublin 4 \\ Ireland}
%\footnote{Research supported by the Claude
%Shannon Institute, Science Foundation Ireland Grant 06/MI/006}\\
%\and Axel Kohnert \\ \inst{Mathematics Department \\ University of Bayreuth \\
%D-95440 Bayreuth\\ Germany}\\ 
%\and Vitaly Skachek  \\ \inst{School of Mathematical Sciences\\ 
%University College Dublin \\ Belfield, Dublin 4 \\ Ireland}
%\footnote{Research supported by the Claude
%Shannon Institute, Science Foundation Ireland Grant 06/MI/006}
%}
\date{}
\maketitle

%\inst{4}}
%\institute{
%School of Mathematical Sciences and Claude Shannon Institute,\\ 
%University College Dublin, Belfield, Dublin 4, Ireland.
%\email{ebyrne@ucd.ie}
%\and 
%School of Mathematical Sciences and Claude Shannon Institute,\\ 
%University College Dublin, Belfield, Dublin 4, Ireland.
%\email{marcus.greferath@ucd.ie}
 %\\
 %\email{ebyrne@ucd.ie,marcus.greferath@ucd.ie,vitaly.skachek@ucd.ie}
%\and
%Mathematics Department, University of Bayreuth, \\
%D-95440 Bayreuth, Germany.
%\email{axel.kohnert@uni-bayreuth.de}
%\and
%School of Mathematical Sciences and Claude Shannon Institute,\\ 
%University College Dublin, Belfield, Dublin 4, Ireland.
%\email{%ebyrne@ucd.ie,marcus.greferath@ucd.ie,
%axel.kohnert@uni-bayreuth.de
%vitaly.skachek@ucd.ie}
%}
%}

%\date{}

%\maketitle

%\begin{document}
%\maketitle

\begin{abstract}
We give further results on the question of code optimality for linear codes over finite 
Frobenius rings for the homogeneous weight.
This article improves on the existing Plotkin bound derived in an earlier paper \cite{plotkin}, and suggests a version of a Singleton bound. We also present
some families of codes meeting these new bounds.
\end{abstract}

\vspace{0.3cm} {\bf Key Words:} codes over rings, finite Frobenius rings,  
homogeneous weights, Plotkin and Singleton bounds.

\section*{Introduction} 

In the early 1990s interest in algebraic codes over finite rings 
was vastly increased due to the discovery that certain non-linear binary 
codes have ${\mathbb Z}_4$-linear representations
(cf.~\cite{hammons94,nech91}). 
Many papers on the topic have been published since then.
A new weight function called the {\em homogeneous weight\/} was discovered by Heise and Constantinescu \cite{cons95,consheis97} and has since proven to be useful in the context of codes over finite rings. 
Examples of homogeneous weights include the Hamming weight on finite fields
and the Lee weight on ${\mathbb Z}_4$. The homogeneous weight may be viewed as a natural generalisation of the Hamming weight for codes over finite rings.

As in traditional algebraic coding theory, a natural
question when dealing with codes over ring alphabets concerns the
criteria that best measure the quality and determine optimality of a code.  
For this reason, the theory requires the establishment of fundamental bounds
relating the standard parameters of code length, size, minumum distance. 
Many of the classical bounds for codes over finite fields have found
an equivalent expression for finite ring codes for the homogeneous weight.
For example, Plotkin and Elias bounds were given in \cite{plotkin} and 
constructions of Plotkin-optimal codes can be read in \cite{butson}.
In \cite{BGS07}, a linear programming bound was derived. 

In this note we present further bounds for linear codes over finite Frobenius rings for the homogeneous weight. We give a refinement of the Plotkin bound given in \cite{plotkin}. 
We also suggest a Singleton-like bound.

\section{Technical Preliminaries} \label{s:prep}

In all that follows, let $R$ be a finite ring with identity. 
The character group of the additive group of $R$ is denoted by 
$\widehat{R}:={\rm Hom}_{\mathbb Z}(R,{\mathbb C}^\times)$. This group
has the structure of an $R$-$R$-bimodule by defining
$\chi^r(x):=\chi(rx)$ and $^r\chi(x):=\chi(xr)$ for all $r,x\in R$,
and for all $\chi\in \widehat{R}$. Summarizing elements from
\cite{wood97a} we come to the following definition:

\begin{definition} \label{def-f-ring}
A finite ring $R$ is called a {\em Frobenius ring\/} if $_R\widehat{R}\cong
{_RR}$.
\end{definition}

It can be seen (cf.~\cite{wood97a}) that if $R$ is a finite Frobenius ring, then $R$ and
$\widehat{R}$ are isomorphic also as right $R$-modules. Hence, there exist 
characters $\chi$ and $\psi$ such that $$\widehat{R} \; = \; \{^r\chi\mid r \in R\} \; = \;
\{\psi^r \mid r \in R\}.$$ 
Such characters are called {\em left generating\/} or {\em right generating\/},
respectively.  Moreover, every left generating
character is at the same time right generating, and a character is (left and/or
right) generating if and only if its kernel  does not contain any non-zero left
or right ideal of $R$. 

The class of finite Frobenius rings is quite large, as the following proposition shows. For a proof see \cite{wood97a} and also \cite{bruce}.

\begin{proposition}
\begin{abc}
\item Any finite principal ideal ring is Frobenius.
\item If $R$ and $S$ are Frobenius ring, then so is $R\times S$.
\item If $R$ is a Frobenius ring, then so is $M_n(R)$, the ring of all $n\times n$-matrices over $R$.
\item If $R$ is a Frobenius ring, and $G$ a finite group, then the group ring $R[G]$ is again a Frobenius ring.
\end{abc}
\end{proposition}

\subsection*{Weight Functions}

The Hamming weight of a word $c \in R^n$ counts the number of the nonzero components of $c$, and hence gives the size of $\supp(c)$. In a way, it could be considered as the actual length of $c$, and hence, we will denote it by $\ell(c)$. For a code $C\leq {_RR^n}$, we write $\ell(C) := |\supp(C)|$. 

We are aware that this notation deviates from the literature, however we ask the reader to accept it, as it will help to avoid confusion with the homogeneous weight and minimum distance that we are going to present now. %We write $\mathbb R^+$ to denote the non-negative real numbers.

\begin{definition}\label{defhomogen}
  A weight function $w: R \longrightarrow {\mathbb R}$ is
  called \emph{(left) homogeneous}, if $w(0)=0$ and the following is
  true:
  \begin{abc}
  \item[\bf (H1)] If $Rx=Ry$ then $w(x)= w(y)$ for all $x,y\in R$.
  \item[\bf (H2)] There exists a real number $\gamma$ such that
    \begin{equation*}
  \sum_{y\in Rx}w(y) \; =\; \gamma \, |Rx|\qquad\text{for all $x\in
    R\setminus \{0\}$}.
\end{equation*}
\end{abc}
\end{definition}

Homogeneous weights were first introduced by Heise and Constantinescu in \cite{consheis97} for integer residue rings, and later generalised to Frobenius rings in \cite{weighted}, and to arbitrary finite rings in \cite{bruce}.

The number $\gamma$ may be thought of as the {\em average value}
of $w$, and condition {\bf (H2)} simply states that this average is the same on all nonzero principal left ideals.
 
 It was shown in ~\cite[Theorem~1.3]{bruce} that, up to the choice of $\gamma$, every finite ring admits a unique (left) homogeneous weight . Moreover, Honold observed in \cite{honold01} that, provided $R$ is Frobenius, the homogeneous weight will allow for an expression in terms of a generating character.
We let $R^{\times}$ denote the group of units of $R$.
 
\begin{proposition}\label{char-homogen} 
  Let $R$ be a finite Frobenius
  ring with generating character $\chi$. 
  Then the (left) homogeneous weights on $R$ are precisely the functions
  \begin{equation*}
w: R\longrightarrow {\mathbb R}, \quad x \mapsto
  \gamma\Big[1-\frac{1}{|R^{\times}|}\sum_{u\in R^{\times}}
  \chi(xu)\Big]
\end{equation*}
where $\gamma$ is a real number.
\end{proposition}

As an immediate consequence, if $R$ is a finite
Frobenius ring, then every left homogeneous weight is also right
homogeneous with the same average value $\gamma$, since
$$\sum_{u\in R^{\times}}\chi(xu)\; = \; \sum_{u\in R^{\times}}\chi(ux).$$

As we will restrict to Frobenius rings in the sequel we will not distinguish between left and right homogeneous weights any more, and simply refer to homogeneous weights instead. Before we continue, we will give examples of homogeneous weights on various instances of finite Frobenius rings.

\begin{example}\label{exhw}
  \begin{abc}
  \item On every finite field $\F_q$ the Hamming weight
    is a homogeneous weight of average value $\gamma=\frac{q-1}{q}$.
  \item On $\Z_4$ the Lee weight is homogeneous with
    $\gamma=1$.
  \item On a local Frobenius ring $R$ with $q$-element residue field the
    weight
    \begin{equation*}
      w\colon R\longrightarrow {\mathbb R},\quad x \mapsto
      \left\{\begin{array}{ccl}
        0&:& x=0,\\
        \frac{q}{q-1}&:&x\in\soc(R),\; x\neq 0,\\
        1&:& \mbox{otherwise},
      \end{array}\right.
    \end{equation*}
is a homogeneous weight of average value $\gamma=1$.
\item On the ring $R$ of $2 \times 2$ matrices over GF$(2)$ the weight
      \begin{equation*}
      w\colon R\longrightarrow {\mathbb R},\quad x \mapsto
      \left\{\begin{array}{ccl}
        0&:& x=0,\\
        2&:&\mbox{$x$ singular},\; x\neq 0,\\
        1&:& \mbox{otherwise},
      \end{array}\right.
    \end{equation*}
    is a homogeneous weight of average value $\gamma = \frac{3}{2}$.
\end{abc}
\end{example}

As is common in coding theory, a weight $w$ on a finite ring $R$ is additively extended to a weight on the $R$-module ${}_RR^n$, i.e.
$$w(c) \; := \; \sum_{i=1}^n w(c_i),\quad \mbox{for $c\in R^n$}.$$
The minimum weight of a linear code is the minimum non-zero weight of any codeword.   
A linear code of length $n$ and minimum homogeneous weight $d$ will frequently be referred to as an $[n,d]$-code. If $R$ is a finite field then the notion of dimension of a linear code is well defined and we write $[n,k,d]$ to denote
a linear code of length $n$, dimension $k$ and minimum weight $d$. We write $(n,M,d)$ to denote a not necessarily linear code over a finite field of length $n$ and minimum distance $d$ with $M$ words.

\section{Shortened and Residual Codes}

We construct new codes from a given code by shortening and puncturing.
The results of this section will be applied in later sections to derive further bounds.

\begin{lemma}\label{lemcoset}
    Let $C \leq {_RR^n}$ be a linear code, and let $x \in R^n$. Then
    $$ \frac{1}{|C|} \sum_{c \in C} w(x+c) = \gamma \ell(C) + \sum_{i\not\in \supp(C)} w(x_i).$$    
\end{lemma}
\pf We compute    \begin{eqnarray*}
           \frac{1}{|C|}\sum_{c \in C} w(x+c) & = & \frac{1}{|C|} \sum_{c \in C}\sum_{i=1}^n w(x_i + c_i)\\
	   &=& \frac{1}{|C|} \sum_{i=1}^n  \sum_{c \in C}\gamma 
                                                   \Big[ 1 - \frac{1}{|R^\times|} \sum_{u \in R^\times} \chi((x_i+c_i)u)    \Big] \\
           & = & \gamma n -  \gamma \frac{1}{|R^\times|}  \sum_{i=1}^n \sum_{u \in R^\times}  \chi(x_iu) 
                      \frac{1}{|C|}\sum_{c \in C} \chi(c_iu).
\end{eqnarray*}
Clearly the projection of $C$ onto some $i$th coordinate is an ideal of $R$, and
since $\chi$ is a generating character we have 
$$\frac{1}{|C|}\sum_{c\in C}\chi(c_iu)\;=\;\left\{\begin{array}{ccl}
0&: & i \in \supp(C),\\
1 &:& \mbox{otherwise,}
\end{array}\right.$$ 
and hence

\begin{eqnarray*} \frac{1}{|C|}\sum_{c \in C} w(x+c) & = &
  \gamma n -  \gamma \frac{1}{|R^\times|} \sum_{i=1}^n \sum_{u \in R^\times} \chi(x_iu)\\
           & = & \gamma \ell(C)  + \sum_{i \notin \supp(C)}  \gamma\Big[1-\frac{1}{|R^\times|} 
           \sum_{u \in R^\times} \chi(x_iu) \Big]\\
           & = & \gamma \ell(C) +\sum_{i\not\in\supp(C)} w(x_i),                                      
     \end{eqnarray*}
     which was the claim.\qed

Given a linear code $C\leq {_RR^n}$ and a subset $S\subseteq \{1,\ldots,n\}$, we define the code 
$$\sh(C,S)\; :=\; \{ c \in C \mid\supp (c) \subset S \},$$ 
which is essentially (namely up to omitting vanishing coordinates) a shortened code. Moreover, we define the residual code 
$$\res(C,S)\; := \{ (c_i)_{i \notin S} \mid c \in C \}.$$

Denoting by $\pi_S$ the projection of $R^n$ onto the coordinates {\em not\/} contained in $S$, it is clear that $\sh(C,S) = \ker(\pi_S) \cap C$ and $\res(C,S) = \pi_S(C)$. Obviously, these codes are related by $C / \sh(C,S) \cong \res(C,S)$.

Finally, for arbitrary $x\in R^n$, for the sake of simplicity of notation we write $\sh(C,x)$ to mean 
$\sh(C,\supp(x))$ and $\res(C,x)$ in place of $\res(C,\supp(x))$. Likewise we will write $\pi_x$ where 
$\pi_{\supp(x)}$ is meant.

In general there is no relationship between $\sh(C,x)$ and $Rx$, except that for $x\in C$ there holds $\sh(C,x)\geq Rx$. The following lemma gives a condition for equality in this containment.

\begin{lemma} \label{lemsho}
     Let $C\leq {_RR^n}$ be a linear code of homogeneous minimum weight $d$, and let $c$ be a word in $C$ that  satisfies $\gamma \ell(c) < d$. Then $\sh(C,c) = Rc$.
    
\end{lemma}
\pf Assuming that there exists $x \in \sh(C,c)$ that is not contained in $Rc$ we first observe that $0 \not\in x+Rc$, which implies $$ d\; \leq\; \frac{1}{|Rc|} \sum_{y\in Rc} w(x+y).$$ We have
     $\supp (x) \subseteq \supp (c)$ and thus may use Lemma \ref{lemcoset} to observe
     $$ \frac{1}{|Rc|} \sum_{y\in Rc} w(x+y) = \gamma \ell (c) + \sum_{i\not\in\supp(c)} w(x_i) = \gamma \ell(c) < d,$$
     which is a contradiction showing the claim.\qed

\begin{corollary} \label{propres}
     Let $C\leq {_RR^n}$ be of homogeneous minimum weight $d$, and let $c \in C$ satisfy $\gamma \ell(c) < d$. Then
     $\res(C,c)$ is of length $n-\ell(c)$, homogeneous minimum weight at least $ d-\gamma\ell(c)$, and satisfies
     $|\res(C,c)|=|C|/|Rc|$.
\end{corollary}

\pf  Let $\res(C,c)$ have minimum homogeneous weight $d'$, and let $x \in C$ such that $w(\pi_c(x))$ assumes $d'$. 
     Then, invoking Lemma \ref{lemcoset}, we have
     $$d \; \leq \; \gamma \ell(c) + \sum_{i\not\in\supp(c)}w(x_i)  \;= \;\gamma \ell(c) + w(\pi_c(x))  \;= \;\gamma \ell(c) + d', $$
     which yields $d' \geq d - \gamma \ell(c)$. Our claim regarding the size of $\res(C,c)$ follows from the fact that $\sh(C,c)= Rc$.\qed

\begin{example}
     Let $C$ be the ${\mathbb Z}_4$-linear Octacode generated by
     $$\left[  \begin{array}{cccccccc}
                                           1 & 0 & 0 & 0 & 3 & 1 & 2 & 1\\
                                           0 & 1 & 0 & 0 & 1 & 2 & 3 & 1\\
                                           0 & 0 & 1 & 0 & 3 & 3 & 3 & 2\\
                                           0 & 0 & 0 & 1 & 2 & 3 & 1 & 1
                    \end{array}    
         \right].$$
         The code $C$ has $256$ words and minimum Lee distance 6 (cf.~\cite{hammons94}). 
         It contains the word $c=[0,0,0,2,0,2,2,2]$ which satisfies $\gamma\ell(c)=4 < 6=d$ where we recall that the Lee weight is homogeneous with  $\gamma=1$.
         Clearly, $|Rc| = 2$ and we puncture
 $C$ on the coordinates $4,6,7,8$ to obtain $\res(C,c)$, which by 
         Corollary \ref{propres} is a linear $[4,d'\geq 2]$ code of size 128. 
 Considering the Gray image (cf.~\cite{hammons94}) of $\res(C,c)$ we arrive at an $(8,128,\geq 2)$ code that obviously meets the (traditional) Singleton bound. This shows that $d'=2$ and hence, $\res(C,c)$ is an optimal code.  
\end{example}

\section{A Refinement of the Plotkin Bound}

  If a linear code $C\leq {_RR^n}$ has maximal support, meaning $\ell(C)=n$, then by observations in \cite{plotkin} or by applying Lemma \ref{lemcoset} we find
      \begin{equation}\label{eqav}
      \frac{|C| -1}{|C|}\, d \; \leq \; \frac{1}{|C|} \sum_{c \in C} w(c)\; =\; \gamma n.
      \end{equation}
            
     We combine this observation with the following theorem to obtain a Plotkin-like bound
     for linear codes.       
     
     \begin{theorem}\label{thplotkin}
             Let $C\leq {_RR^n}$ be a linear $[n,d]$ code satisfying 
             $\gamma n < d$, and let $c\in C$ be such that $\gamma\ell(c) < d$.
             Then there holds
             $$ |C| \; \leq 
	    \;  |Rc|\, \frac{d-\gamma \ell(c)}{d - \gamma n}.$$      \end{theorem}
      
\pf   Suppose that $C_1:=\res(C,c)$ has length $n_1$
      and minimum homogeneous weight $d_1$.
      From (\ref{eqav}) and Corollary \ref{propres} we have
      \begin{equation*}
            n \; = \; \ell(c) +  n_1 \; \geq \;  \ell(c)+  \frac{|C_1|-1}{|C_1|}  \,\frac{d_1}{\gamma} 
                \; \geq \; \ell(c) +  \frac{|C_1|-1}{|C_1|} \Big( \frac{d}{\gamma}- \ell(c) \Big) \\
      \end{equation*}
      From Corollary \ref{propres} we know that $|C_1| = |C|/|Rc|$, which gives   
      \begin{eqnarray*}        
                n & \geq & \ell(c) +  \Big( 1 - \frac{|Rc|}{|C|} \Big)\Big( \frac{d}{\gamma}- \ell(c) \Big).\\
               %& \geq & l \frac{|R|}{|C|} + \Big[ 1- \frac{|R|}{|C|} \Big] \frac{d}{\gamma}
      \end{eqnarray*}
      Rearranging this inequality yields the result.      
      \qed

      \begin{example}
      Let $m\in {\mathbb N}$ and let $n = m\times (|R|^m-1)$. We consider the code $C \leq ~_RR^n$ 
      which is generated by the 
      $m\times n$ matrix $G$ whose columns comprise the distinct nonzero elements of $R^m$.
      It is not difficult to see that $C$ is a constant weight code of homogeneous weight $\gamma |R|^m$. If namely $x \in R^m$ then
      \begin{eqnarray*}
         w(xG) & = & \sum_{g \in R^m} w(x \cdot g)
                =  \sum_{g \in R^m} \gamma \Big[ 1 - \frac{1}{|R^\times|}\sum_{u \in R^\times} \chi(ux\cdot g) \Big]\\
               & = & \gamma \Big[ |R|^m - \frac{1}{|R^\times|}\sum_{u \in R^\times}\sum_{g \in R^m}\chi(ux\cdot g) \Big] \\
               & = & \left\{ \begin{array}{ccl}
                                        0        &: & x = 0 \\
                                        \gamma |R|^m &: & \mbox{otherwise.}
                                   \end{array} 
                          \right.
      \end{eqnarray*}
     Moreover, $n = |R|^m-1 < |R|^m = \frac{d}{\gamma}$. It can also be shown that 
   $\ell(c) \leq n < \frac{d}{\gamma}$ for each word $c\in C$.
             The Hamming weight
             of an arbitrary word $c=xG$ of $C$ corresponds to the size of the annhilator submodule
             $x^\perp =  \{ y \in R^m \mid x \cdot y = 0 \} \leq R^m_R$ by the equation
             $$\ell(c) \; = \; |R|^m - |x^\perp| \; = \; |R|^m - \frac{|R|^m}{|Rc|}.$$
             Therefore, the upper bound on $|C|$ determined by Theorem \ref{thplotkin} is
             $$ |C| \; \leq \;  |Rc|\,  \frac{d-\gamma \ell(c)}{d - \gamma n} \;  = \;             |Rc| \Big[ |R|^m -|R|^m +  \frac{|R|^m}{|Rc|} \Big]\; = \; |R|^m,$$
             which is met sharply by $C$.
      \end{example}
      
      We will refer to the code in the preceding example as a {\em Simplex code\/}.

      \begin{corollary}\label{corplotkin} %eimo why possibly equalities?
            Let $C \leq {_RR^n}$ be of minimum homogeneous weight $d$ and minimum Hamming weight $\ell$ where $\ell \leq n \leq \frac{d}{\gamma}$. Then
            $$ |C| \; \leq \; |R| \,\frac{d-\gamma \ell}{d - \gamma n}.$$
      \end{corollary}
      
      It is straightforward to verify that for linear codes, this gives a refinement of the Plotkin bound given in 
      \cite{plotkin} for $\ell < \frac{d}{\gamma} < \ell \frac{|R|}{|R|-1}$.

      In fact we can do even better, taking into account some properties of $R$. For this, we first make an elementary but useful observation.

      \begin{lemma}\label{lemminham}
      Let $C \leq {_R}R^n$ and let $c\in C$ have minimum Hamming weight in $C$. 
      Then there exists $\alpha\in R$ and a family $(u_i)_{i\in\supp(c)}$ of invertible elements of $R$ such that 
      $c_i = \alpha u_i$ for all $i\in\supp(c)$. 
      In particular, $Rc \cong R\alpha$. 
      \end{lemma}
      
\pf
      Since $c$ is of minimal Hamming weight, we have
      $\ell(\lambda c) = \ell(c) $ for each $\lambda \in R$, unless $\lambda c=0$. For this reason, the left annihilators $c_i^\perp:=\{\lambda \in R \mid \lambda c_i=0\}$ must all be the same for $i\in \supp(c)$, which holds if and only if 
      the $c_i R$ coincide for all such $i\in {\rm supp}(c)$. Then the claim follows from \cite[Thm 5.1]{wood97a}.\qed
 
 \begin{lemma}\label{marcus}
      Let $C\leq {_RR^n}$ be a linear code of minimum homogeneous weight $d$ and minimum Hamming weight $\ell$ where $\gamma\ell < d$. If $c\in C$ is a word of minimum Hamming weight then $Rc$ is a simple submodule of $C$.
      \end{lemma}
\pf Suppose that $Rc' \leq Rc$ for some nonzero $c' \in C$. Then $\ell(c)=\ell(c')$ and in particular $\supp( c') = \supp( c)$.
      By Lemma \ref{lemsho}, we find that $Rc'= \sh(C,c') = \sh(C,c) = Rc$.
      Thus $Rc$ is a simple submodule of $C$.\qed
   
      \begin{corollary}\label{corminham}
      Let $C\leq {_RR^n}$ be a linear code of minimum homogeneous weight $d$ and minimum Hamming weight $\ell$ where $\ell < n \leq \frac{d}{\gamma}$. Let $Q$ be the maximum size of any minimal ideal
      of $R$. Then
      $$ |C| \leq Q \,\frac{d-\gamma \ell}{d - \gamma n}.$$
      \end{corollary}

\pf
      Let $c \in C$ be of Hamming weight $\ell$. By the preceding lemma we know that $Rc$ is a simple submodule of $C$. Combining this with Corollary \ref{corplotkin} the claim follows immediately.\qed

      \begin{example}
      We again study the Simplex Code, this time over the ring $R$ of all $2\times 2$-matrices over ${\mathbb F}_2$. This code is of length $n=16^m-1$ for suitable $m$, and its minimum Hamming weight of is $16^m - \frac{16^m}{4} = \frac{3}{4}16^m$. The ring $R$ has $3$ minimal ideals, each of size $4$, and
           so, from Corollary \ref{corminham}, we have
           $$16^m \;= \;|C| \;\leq \; 4 \,\frac{16^m \gamma - \frac{3}{4}16^m \gamma}{16^m \gamma - (16^m-1) \gamma}  \;=\; 4 \frac{16^m}{4} \;=\; 16^m,$$ showing that the bound in the previous corollary is met sharply.
      \end{example}

\subsection{A Singleton bound}   

       Let $C$ be an $[n,d]$ code over $R$ satisfying $n \leq \frac{d}{\gamma}$.
       If $c \in C$ is a codeword satisfying $\ell:=\ell(c) < n \leq \frac{d}{\gamma}$ then by  Corollary \ref{propres} we see that   
       $C_1:=\res(C,c)$ is an $[n_1,d_1]$ code
       over $R$, isomorphic to $C / Rc$ with $d_1 \geq d - \gamma \ell$ and
       $$n_1 \;= \; n - \ell \; \leq \;  \frac{d}{\gamma} - \ell \; \leq 
      \;  \frac{d_1}{\gamma} .$$
       
       Setting $C_0:=C$, we construct a sequence of 
       $[n_i,d_i]$ codes $C_i$ as follows: 
       for each $i$, as long as there exists $c^i \in C_i$ with Hamming weight $\ell_i:=\ell(c^i)<n_i$, define $C_{i+1}:=\res(C_{i},c^i)$. We observe that $n \leq \frac{d}{\gamma}$ implies
       $\ell_i < n_i =n_{i-1} - \ell_i < \frac{d_i}{\gamma}$ for each $i \geq 1$. Therefore, 
       from Lemma \ref{lemsho}
       we have a finite sequence of codes
       $$C_0=C,\; C_1 \cong C_0 /Rc^0,\;C_2 \cong C_1 /Rc^1,...,\; C_r \cong C_{r-1}/Rc^{r-1} $$
       of length $r+1$ for some nonnegative integer $r$.
       Moreover, for each  $i \in \{1,...,r\}$ we have
       \begin{equation}\label{eqci}
        |C_i|\; =\; \frac{|C_{i-1}|}{|Rc^{i-1}|} \;=\; \frac{|C|}{|Rc^{0}|\cdots |Rc^{i-1}|} \quad \mbox{and } 
       d_i \geq d_{i-1}-\gamma \ell_{i-1} > 0.
       \end{equation}
       Note that the final code $C_r$ has the property that each of its non-zero words has constant 
       Hamming weight $n_r$, so taking any further quotients by $c^r \in C_r$ will result in a code of length zero.
%       We may write $C_r =\sh(C_{r},c^{r})$ for any $c^r \in C_r$. Since it may occur that 
%       $\ell(c^r):=\ell_r=n_r =\frac{d_r}{\gamma}$ we cannot apply Lemma \ref{lemsho} to determine 
%       that $C_r=\sh(C_{r},c^{r})=Rc^r$. 
       Employing a simple counting argument (e.g.~traditional Singleton bound for the Hamming distance) it can be shown that $|C_r| \leq |R|$.
              
       From Equation (\ref{eqci}) we have
       \begin{equation}\label{eqC}
       |C| \;=\;   |Rc^0|\,|Rc^1| \cdots |Rc^{r-1}|\,|C_r|.
       \end{equation}
       
       The existence of such a sequence of $r+1$ codes leads to the following inequality. 
      
       \begin{eqnarray}
              n  =  \sum_{i=0}^{r} \ell_i & \geq & \frac{|Rc|-1}{|Rc|} \frac{d}{\gamma} + \sum_{i=1}^{r} \ell_i\\
                                       & \geq & \frac{|Rc|-1}{|Rc|} \frac{d}{\gamma} + r,\label{eqsing}
       \end{eqnarray}
       
       This will yield a type of Singleton bound for the homogeneous weight. First we need one further observation.
       
       \begin{lemma}\label{lemQ}
            Let $C$ be an $[n,d]$ code over $R$ satisfying $\gamma n \leq d$.
            Let $Q:=\max \{ |Rc|\mid c \in C\}$ and let 
            $P:=\max \{ |Rc|\mid  c \in C,\ell(c)< n\}$. 
           If $c \in C$ satisfies $\ell(c) < n$ then $|Rc'| \leq Q$ for each
            $c' \in \res(C,c)$. Moreover, if $\ell(c')<n-\ell(c)$ then $|Rc'|\leq P$.
       \end{lemma}
       
       \pf
        Let $c' \in \res(C,c)$. From Lemma \ref{lemsho}, we have $\res(C,c) \cong C / Rc$ and hence, there is some $x \in C$ such that $Rc' \cong (Rx+Rc) / Rc$. Consequently, 
        $$|Rc'| \; = \; \frac{|Rx+Rc|}{|Rc|} \; = \; \frac{|Rx|}{|Rx \cap Rc|} \; \leq \; |Rx| \; \leq \; Q.$$
        If $|Rc'|>P$ then $|Rx|>P$ and hence $\ell(x) = n$, which implies $\ell(c')=n-\ell(c)$.\qed
       
       \begin{theorem}\label{thsing} %obacht
          Let $C$ be an $[n,d]$ code over $R$ satisfying $\gamma n \leq d$ and with minimum Hamming weight less than $n$.
          Let $P:=\max \{ |Rc|\mid c \in C,\ell(c)<n\}$. 
          %$$Q:= \left\{ \begin{array}{cll}
           %                      |R| &, & \mbox{if }$\ell(C) =n$\\%\mbox{ has constant Hamming weight }$n$\\
            %                     \max \{|Ra|: a \in R^n, Ra \leq C, \ell(a) < n \}&, & {\mbox{otherwise}} .
            %                        \end{array}
            %                        \right.$$
          Then 
          $$ n - \left\lceil \frac{P-1}{P} \frac{d}{\gamma} \right\rceil \;\geq \;\left\lceil \log_{P} |C| - \log_{P} |R| \right\rceil.$$
          %Moreover, if $\frac{d}{\gamma} < n$ then
          %$$ n - \left\lceil \frac{Q-1}{Q} \frac{d}{\gamma} \right\rceil \geq \left\lceil \log_{Q} |C| - 1 \right\rceil.$$
       \end{theorem}

        \pf
        Let $c \in C$ such that $|Rc|=P$. % such that $\ell(c) < n$. 
        With the same notation as before, from Lemma \ref{lemsho} and Corollary \ref{propres}, 
        there exists a sequence of words
        $c=c^0,c^1,...,c^{r-1}$ and linear codes $C=C_0,C_1,...,C_r$ such that, for $i=1,...,r$,
        $C_i:=\res(C_{i-1},c^{i-1})$ is an $[n_i,d_i]$ code, and for $i=0,...,r-1$,
        $c^i \in C_i$, $\ell(c^i) < n_i \leq \frac{d_i}{\gamma}$ 
        and 
        $C_{i} \cong C_{i-1} / Rc^{i-1}$.
        %For each $i \in \{1,...,r-1\}$, there is some $v^i \in C_{i-1}$ such that
        %$Rc^i \cong (Rv^{i}+Rc^{i-1}) / Rc^{i-1} $, %If $|Rc'| \leq Q$ 
        %and hence 
        %$$|Rc^i| = \frac{|Rv^i+Rc^{i-1}|}{|Rc^{i-1}|}  = \frac{|Rv^i|}{|Rv^i \cap Rc^{i-1}|} \leq |Rv^i| \leq Q.$$
        %If $|Rc'| > Q$, then $\ell(c')=n$, and hence $\ell(c_1) =n_1$, in which case $C_1$ thus
        %$l(c_1) < n_1$ implies $|Rc_1| \leq Q.$
        As observed in Lemma \ref{lemQ}, we have $|Rc^i| \leq P$ for $i=1,...,r-1$.         The code $C_r = \sh(C_r,c^r)$ has constant Hamming weight $n_r$ and hence $|C_r|\leq |R|$. 
        Then
        $$|C| \; = \;  |Rc|\,|Rc^1|\cdots |Rc^{r-1}|\,|C_r| \; \leq P^{r}\,|R|,$$               so clearly $r \geq \lceil \log_{P} |C| -\log_P |R| \rceil.$ 
        The inequality in (\ref{eqsing}) gives
        $$ n - \left\lceil \frac{P-1}{P} \frac{d}{\gamma} \right\rceil \;\geq\; \left\lceil \log_{P} |C| - \log_{P} |R| \right\rceil.$$\qed
        
        \begin{corollary}
        Let $C$ be an $[n,d]$ code over $R$ satisfying $n <  \frac{d}{\gamma}$,
        and let $Q:=\max \{|Rc| \mid c \in C\}$. Then
        $$ n - \left\lceil \frac{Q-1}{Q} \frac{d}{\gamma} \right\rceil \;\geq\; \left\lceil \log_{Q} |C| -1 \right\rceil. $$      
        \end{corollary}
                \pf 
        Let $c \in C$ such that $|Rc|=Q$. As before, we recursively define a sequence of 
        $[n_i,d_i]$ codes
        $C_i:=\res(C_{i-1},c^{i-1})$ with $C_1:=\res(C,c)$,
        $c^i \in C_i$, $\ell(c^i) < n_i \leq \frac{d_i}{\gamma}$ 
        and $C_{i} \cong C_{i-1} / Rc^{i-1}$.
        Now $n < \frac{d}{\gamma}$ implies $n_r < \frac{d_r}{\gamma}$ 
        so from Lemma \ref{lemsho} we have $C_r =\sh(C_r,c^r)=Rc^r$. Then 
        $|C_r| \leq Q$ and hence $|C|\leq Q^{r+1}$. Then $r \geq \left\lceil \log_{Q} |C| - 1 \right\rceil$
        and again the result follows from the inequality in (\ref{eqsing}).\qed  
        
        %If $C$ has constant Hamming weight $n$ then $|C| \leq |R|$ and so the result follows immediately from (\ref{eqav}).

        We may deduce the following weaker result directly from Equation (\ref{eqC}).

        \begin{proposition}
            Let $C \leq ~ _RR^n$ be an $[n,d]$ linear code and suppose that $\gamma n \leq d$.
            Then 
            $$ n - \left\lceil \frac{|R|-1}{|R|} \frac{d}{\gamma} \right\rceil \; \geq \; \left\lceil \log_{|R|} |C| - 1 \right\rceil.$$
        \end{proposition}
      
      We give an example of what could be called an MDS code over a finite chain ring $R$, using points from a {\em projective Hjelmslev geometry}.
      
      \begin{example}
            Let $R$ be a chain ring of length 2 with $q$-element residual field.
            Then $R^\times = R \backslash \rad (R)$ and $|R| = q^2$. Let $F:=R^2 \backslash \rad (R^2)$.
            We denote by  ${\rm PHG}(R^2)$ the projective Hjelmslev line with point set ${\cal P} := \{ xR \mid x \in F \}$.
            Note that ${\cal P}$ contains $q^2+q$ distinct points (cf. \cite[p. 83]{geom}). 
           
	   For $n:=q^2+q$ let $C \leq {_R}R^n$ be the code generated by the $2 \times n$ generator matrix 
            $G = [g_1,...,g_n]$ whose
            columns comprise elements of $R^2$ corresponding to distinct points in ${\cal P}$.
            Clearly $\ell(c) <n$ for each $c \in C$. Moreover,
            $C$ is free of rank $2$ and the maximal cyclic submodules of $C$ have size $P:=|R|=q^2$. 
            With $r = \lceil \log_{P} |C| -1 \rceil = \log_{q^2} q^4 -1 = 1$ and $\gamma=1$, 
            each word $xG$ of $C$ has weight 
            $$w(xG) =  |J_1| + \frac{q}{q-1}|J_2| = \left\{ \begin{array}{lcl}
                                                                            q^2 + \frac{q}{q-1}(q-1) = q^2 + q &:& x \in F  \\ 
                                                                            q^2 \frac{q}{q-1} = \frac{q^3}{q-1} &:& x \in \rad (R^2), \; x\neq 0 
                                                                            \end{array} \right.,$$
            where $J_1=\{ j \mid x \cdot g_j \in R^\times \}$ and 
            $J_2=\{ j \mid x \cdot y_j \in \rad(R) \setminus \{ 0\} \}$.
            Then $d=n=q^2+q$ and
            \begin{eqnarray*}
                                    n -   \left\lceil \frac{q^2-1}{q^2}\, d \right\rceil 
                                    & = & n - \left\lceil \frac{q^2-1}{q^2}\, (q^2+q) \right\rceil 
                                     =   n - \left\lceil q^2 + q - 1 - \frac{1}{q} \right\rceil \\
                                    & = & q^2 + q - q^2 - q + 1 = 1 = r,
           \end{eqnarray*}                         
             which meets the bound given in Theorem \ref{thsing}.                       
      \end{example}

%%%%%%%%%%%%%%%%%%%%%%%%%%%%%%%%%%%%%%%%   

\renewcommand{\baselinestretch}{0.95}
\normalsize
       
\providecommand{\bysame}{\leavevmode\hbox to3em{\hrulefill}\thinspace}


\begin{thebibliography}{10}



\bibitem{BGS07}
E. Byrne, M.~Greferath and M.~E.~O'Sullivan, {\em The linear programming bound for codes 
over finite Frobenius rings,} Designs, Codes and Cryptography,
Vol. 42 ,  {\textbf 3}  (2007), pp. 289 - 301.

\bibitem{cons95}
I.~Constantinescu, \emph{Lineare {C}odes \"uber {R}estklassenringen ganzer
{Z}ahlen und ihre {A}utomorphismen bez\"uglich einer verallgemeinerten
{H}amming-{M}etrik}, Ph.D. thesis, Technische Universit\"at M\"unchen, 1995.



\bibitem{consheis97}
I.~Constantinescu and W.~Heise, \emph{A metric for codes over residue class
rings of integers}, Problemy Peredachi Informatsii \textbf{33} (1997), no.~3,
22--28.


\bibitem{bruce}
M.~Greferath and S.~E.~Schmidt, \emph{Finite-Ring Combinatorics and MacWilliams Equivalence Theorem}, J. of Combinatorial Theory (A) \textbf{92} (2000), 17--28.

\bibitem{butson}
M. Greferath, G. McGuire, M. E. O'Sullivan, {\em On Plotkin Optimal Codes over Finite Frobenius Rings}, {\em Journal of Algebra and Its Applications} \textbf{5} (2006), no. 6, 799--815. 

\bibitem{plotkin}
M.~Greferath and M.~E.~O'Sullivan, \emph{On Bounds for Codes over Frobenius Rings under Homogeneous Weights}, Discrete Mathematics \textbf{289} (2004), 11--24.

\bibitem{hammons94}
{A. R.} Hammons, {P. V.} Kumar, {A. R.} Calderbank, {N. J. A. } Sloane, and
{P.} Sol\'e, \emph{The ${\mathbb Z}_4$-linearity of {K}erdock, {P}reparata,
{G}oethals and related codes}, IEEE Trans. Inform. Theory \textbf{40} (1994),
301--319.

\bibitem{weighted}
{W.} Heise, {T.} Honold, {A. A.} Nechaev, \emph{Weighted modules and
  representations of codes}, Proceedings of the ACCT 6 (Pskov, Russia, 1998),
   123-129.

\bibitem{honold01}
{T.} Honold, \emph{A characterization of finite Frobenius rings}, Arch. Math. (Basel), \textbf{76} (2001), 406--415.

\bibitem{geom}
R. Kaya, P. Plaumann, K. Strambach, {\em Rings and Geometry}, NATO ASI Series, Reidel, (1984).

\bibitem{nech91}
{A. A.} Nechaev, \emph{Kerdock codes in a cyclic form}, Discrete Math. Appl. \textbf{1} (1991), 365--384.

\bibitem{wood97a}
{J. A.} Wood, \emph{Duality for modules over finite rings and applications to
coding theory}, Amer. J. Math. \textbf{121} (1999), 555--575. 

%\bibitem{veld} 
%G.~T\"orner and F.~D.~Veldkamp, \emph{Literature on geometry over rings.} J. Geom.  {\bf 42}  (1991),  no. 1-2, 180--200.


\end{thebibliography}
\end{document}